\documentclass[centertags,12pt]{amsart}
\usepackage{latexsym}
\usepackage{amsthm}
\usepackage{amssymb}
\usepackage{color}
\usepackage[utf8]{inputenc}
\usepackage{mathrsfs}
\usepackage[english]{babel}
\usepackage{setspace}
\usepackage{ textcomp }
\usepackage{graphicx}
\usepackage{float}
\usepackage{tikz}
\usepackage{soul}
\graphicspath{ {images/} }

\numberwithin{equation}{section}
\makeatletter
\renewcommand{\subsection}{\@startsection
	{subsection}{2}{0mm}{\baselineskip}{-0.25cm}
	{\normalfont\normalsize\bf}}
\makeatother

\newtheorem{theorem}{Theorem}[section]
\newtheorem{proposition}[theorem]{Proposition}
\newtheorem{lemma}[theorem]{Lemma}
\newtheorem{corollary}[theorem]{Corollary}

{\theoremstyle{definition}

	}
\theoremstyle{remark}
\newtheorem{remark}[theorem]{Remark}
\newtheorem{obs}[theorem]{Observation}

\newcommand{\bfq}{{\mathbb F_q}}

\newcommand{\PGU}{{\rm PGU}}

\newcommand{\cX}{{\mathcal X}}

\newcommand{\cH}{{\mathcal H}}

\newcommand{\cD}{{\mathcal D}}

\newcommand{\GK}{{\mathcal{GK}}}

\newcommand{\Qbar}{{\bar{Q}}}
\newcommand{\Pbar}{{\bar{P}}}
\newcommand{\Rbar}{{\bar{R}}}

\begin{document}
	
	\author[Maria Montanucci and Vicenzo Pallozzi Lavorante]{Maria Montanucci and Vicenzo Pallozzi Lavorante}
	\title{AG codes from the second generalization of the GK maximal curve}
	
	\begin{abstract}
		The second generalized $GK$ maximal curves $\GK_{2,n}$ \cite{BM}
		are maximal curves over finite fields with $q^{2n}$ elements, where $q$ is a prime power and $n \geq 3$ an odd integer. In this paper we determine the structure of the Weierstrass semigroup $H(P)$ where $P$ is an arbitrary $\mathbb{F}_{q^2}$-rational point of $\GK_{2,n}$. We show that these points are Weierstrass points and the Frobenius dimension of $\GK_{2,n}$ is computed. A new proof of the fact that the first and the second generalized $GK$ curves are not isomorphic for any $n \geq 5$ is obtained. AG codes and AG quantum codes from the curve $GK_{2,n}$ are constructed; in some cases, they have  better parameters with respect to those already known. 
	\end{abstract}
	
	\maketitle
	
	\begin{small}
		
		
		
	\end{small}
	
	\section{Introduction}
	
	Let $p$ be a prime and $q=p^n$ a prime power. We denote by $\bfq$ the finite field with $q$ elements. In this paper, by algebraic curve we mean a projective, geometrically irreducible, non-singular algebraic curve defined over $\bfq$. We say that an algebraic curve $\cX$ is $\bfq$-maximal if its number of $\bfq$-rational points $|\cX(\bfq)|$ attains the Hasse-Weil upper bound, namely $ |\cX(\bfq)|=q^2+1+2g\sqrt{q} $, where $g$ is the genus of $\cX$. 
	Apart from being mathematical objects with intrinsic interest,  maximal curves are often used for applications in Coding Theory.
	Algebraic-geometric codes (AG codes) are error correcting codes constructed from algebraic curves, introduced by Goppa in the '80s; see \cite{Goppa}. Roughly speaking, AG codes have better parameters when the underlying curve has many rational points. In this context, maximal curves play a key role being the curves having the largest possible number of rational points with respect to their genus.
	
	The most important tool for constructing AG codes is the Weierstrass semigroup $H(P)$ of $\cX$ at $P$. 
	The semigroup $H(P)$ is defined to be the set of all integers $k$ for which there exists a rational function on $\mathcal{X}$ having pole divisor $kP$. 
	Clearly $H(P)$ is a subset of $\mathbb{N}=\{0,1,2,\ldots\}$. The Weierstrass gap Theorem \cite[Theorem 1.6.8]{Sti}, states that the set $G(P):= \mathbb{N} \setminus H(P)$ contains exactly $g$ elements, called \textit{gaps}. 
	In the finite field setting, the parameters of AG codes constructed from $\cX$ rely on the inner structure of the semigroup $H(P)$; see e.g. \cite{TV1991}.
	The structure of $H(P)$ is not always the same for every point $P$ of $\mathcal{X}$. However, the semigroup $H(P)$ is known to be the same for a generic point $P \in \mathcal{X}$, but there may exist a finite number of points on $\mathcal{X}$, called Weierstrass points, with a different set of gaps. These points are of independent interest, for example in St\"ohr-Voloch Theory \cite{SV}, but they are also relevant in the study of AG codes. 
	Indeed,  most of the codes constructed from maximal curves are those having the best parameters known in the literature. Furthermore, maximal curves often have large automorphism groups which in many cases can be inherited by the AG code itself: this can bring good performances in encoding \cite{Joyner2005} and decoding \cite{HLS1995}.
	
	Recently, Beelen and Montanucci in \cite{BM} introduced a new infinite family of maximal curves $\GK_{2,n}$. For any odd $n \geq 3$ the curve $\GK_{2,n}$ is $\mathbb{F}_{q^{2n}}$-maximal and for $n=3$ is isomorphic to the well-known $GK$ maximal curve constructed by Giulietti and Korchm\'aros in \cite{GK}. The $GK$ curve was already generalized by Garcia, G\"uneri and Stichtenoth \cite{GGS} to an infinite family of $\mathbb{F}_{q^{2n}}$-maximal curves $\GK_{1,n}$ where $n \geq 3$ is an odd prime. For this reason we will refer to the maximal curve $\GK_{2,n}$ as the second generalization of the GK maximal curve.
	
	The aim of this paper is to investigate Weierstrass semigroups at some  points of the curve $\GK_{2,n}$ and to use them to construct several examples of AG codes and quantum codes with good parameters. 
	
	More precisely, in Sections \ref{resO1} and \ref{resO2}  the structure of the Weierstrass semigroups $H(P)$ where $P$ is an arbitrary $\mathbb{F}_{q^2}$-rational point of $\GK_{2,n}$ is considered. The following result will be proven; see also Theorems \ref{primo} and \ref{secondo}.
	
	\begin{theorem} 
		Let $q$ be a prime power, $n \geq 5$ odd and $m=(q^n+1)/(q+1)$. If $P \in \GK_{2,n}(\mathbb{F}_{q^2})$ then the Weierstrass semigroup $H(P)$ at $P$ is
		\begin{itemize}
			\item $H(P)=\langle mq+i(q^2-q),q^n+1 | i=0, \dots, s \rangle $, if $P$ is an ideal point of $\GK_{2,n}$;
			\item $H(P)=\langle q^n+1-m,q^n+1-k \mid k=0,\ldots, (m-1)/(q^2-q)\rangle,$ otherwise.
		\end{itemize}
	\end{theorem}
	
	Section \ref{secFrob} is devoted to the computation of an important birational invariant of $\GK_{2,n}$, namely its Frobenius dimension; see Theorem \ref{frobx}. 
	\begin{theorem}
		Let $q$ be a prime power and $n\geq 5$ an odd integer. The Frobenius dimension of the second generalized $GK$ curves $\GK_{2,n}$ is
		\begin{equation*}\label{rdim}
		r=\frac{m-1}{q^2-q}+2, \mbox{ where } m=\frac{q^n+1}{q+1}.
		\end{equation*}
	\end{theorem}
	
	As an application the following corollary will be proven; see also Corollary \ref{wp} and Theorem \ref{isomor}.
	
	\begin{corollary}
		Let $q$ be a prime power and $n \geq 3$ be odd.
		\begin{itemize}
			\item The generalized GK curve $ \GK_{2,n}$ is isomorphic to $\GK_{1,n}$ if and only if $n =3$.
			\item The $\mathbb{F}_{q^{2n}}$-rational points of $\GK_{2,n}$ are Weierstrass points. 
		\end{itemize}
	\end{corollary}

	Finally, in Section \ref{appAG}  we apply our results to the construction of AG codes and AG quantum codes from the curve $\GK_{2,n}$.
	
	Comparisons  with codes already known in the literature will be provided, pointing out how in some cases  better parameters are obtained; see Remark \ref{rem6.3}.

	\section{Preliminary results}
	
	From a result commonly known as Kleiman-Serre covering result \cite{KS}, we know that every curve which is $\mathbb{F}_{q^2}$-covered by an $\mathbb{F}_{q^2}$-maximal curve is itself also $\mathbb{F}_{q^2}$-maximal. The most important example of $\mathbb{F}_{q^2}$-maximal curve is the Hermitian curve $\cH_q$, with affine equation \[Y^{q+1}=X^{q+1}-1.\] The automorphisms group of $\cH_q$ is very large compared to $g(\cH_q)$. Indeed it is isomorphic to $\PGU(3,q)$ and its order is larger than $16 g(\cH_q)^4$. Moreover $\cH_q$ has the largest genus admissible for an $\mathbb{F}_{q^2}$-maximal curve and it is the unique curve having this property up to birational isomorphism, see \cite{RS}.
	Few examples of maximal curves not covered by $\cH_q$ are known in the literature. In \cite{GK} Giulietti and Korchm\'aros constructed an $\mathbb{F}_{q^6}$-maximal curve  which is not a subcover of the Hermitian curve $\cH_{q^3}$. Two generalizations of the Giulietti-Korchmàros curve (GK curve) into infinite families of maximal curves are known in the literature and they are not Galois subcover of the Hermitian curve. 
	The first generalization $\GK_{1,n}$ was introduced by Garcia, G\"{u}neri and Stichtenoth in \cite{GGS}. For  a prime power $q$ and an odd integer $n\geq 3$ the curve $\GK_{1,n}$ is given by the affine space model
	\begin{equation}\label{GGS}
	\GK_{1,n}:\begin{cases}
	Z^m=Y^{q^2}-Y \\
	Y^{q+1}=X^{q}+X
	\end{cases}
	\end{equation} where $m:=\frac{q^n+1}{q+1}$. The curve $\GK_{1,n}$ is $\mathbb{F}_{q^{2n}}$-maximal of genus $g(\GK_{1,n})=(q-1)(q^{n+1}+q^n-q^2)/2$ and $\GK_{1,3}$ is the GK curve.
	
	Recently Beelen and Montanucci \cite{BM} constructed another infinite family $\GK_{2,n}$ of maximal curves generalizing the GK curve. For any prime power $q$ and odd $n\geq3$ the curve $\GK_{2,n}$ is given by \begin{equation}\label{BM}
	\GK_{2,n}:\begin{cases}
	Z^m&=Y\frac{X^{q^2}-X}{X^{q+1}-1} \\
	Y^{q+1}&=X^{q+1}-1
	\end{cases}
	\end{equation} where again $m:=\frac{q^n+1}{q+1}$. The main properties of $\GK_{1,n}$, for a fixed $n$, are summarized in the next propositions.
	
	\begin{proposition} \cite[Section 2]{BM} \label{BMprop}
		Let $\GK_{2,n}$ be defined as above. Then the following holds. \begin{itemize}
			\item[$\bullet$] $\GK_{2,n}$ is an absolutely irreducible $\mathbb{F}_{q^{2n}}$-maximal curve. The genus of $\GK_{2,n}$ is \[g(\GK_{2,n})=\frac{(q-1)(q^{n+1}+q^n-q^2)}{2}.\]  The curve $\GK_{2,3}$ is $\mathbb{F}_{q^6}$-isomorphic to the GK curve. Even though $g(\GK_{1,n})=g(\GK_{2,n})$, the curves $\GK_{1,n}$ and $\GK_{2,n}$ are not isomorphic over $\overline{\mathbb{F}}_{q^{2n}}$ for any $n \geq 5$.
			\item[$\bullet$] For any odd $n \geq 3$ and $q>2$, $\GK_{2,n}$ is not Galois covered by the Hermitian curve $\cH_{q^n}$. If $q=2$ then $\GK_{2,n}$ is Galois covered by $\cH_{q^n}$ over $\mathbb{F}_{q^{2n}}$ for every odd $n \geq 3$.
			\item[$\bullet$] The automorphism group of $\GK_{2,n}$ is isomorphic to $ {\rm SL}(2,q) \rtimes C_{q^n+1}$ where $C_k$ denotes the cyclic group with $k$ elements.
		\end{itemize}
	\end{proposition}
	
	
	Let $\mathbb{F}_{q^{2n}}(x,y,z)$ be the function field of $\GK_{2,n}$ where \[y^{q+1}=x^{q+1}-1 \mbox{ and } z^m=y \dfrac{x^{q^2}-x}{x^{q+1}-1}.\] Then $\mathbb{F}_{q^{2n}}(x,y,z)$ is an extension of the Hermitian function field $\mathbb{F}_{q^{2n}}(x,y)$. The following proposition describes the ramification structure in the function fields extension $\mathbb{F}_{q^{2n}}(x,y,z)/\mathbb{F}_{q^{2n}}(x,y)$ and the short-orbits structure of the automorphism group of $\GK_{2,n}$ over the algebraic closure of $\mathbb{F}_{q^{2n}}$.
	
	\begin{proposition}\cite[Section 4]{BM} \label{ramificati}
		For the function field $\mathbb{F}_{q^{2n}}(x,y,z)$ of $\GK_{2,n}$ the following holds.
		\begin{itemize}
			\item[$\bullet$] The places centered at the $q^3+1$ $\mathbb{F}_{q^2}$-rational points of the Hermitian curve $\cH_q$ are totally ramified in the function fields extension $\mathbb{F}_{q^{2n}}(x,y,z)/\mathbb{F}_{q^{2n}}(x,y)$. Moreover, $\mathbb{F}_{q^{2n}}(x,y,z)/\mathbb{F}_{q^{2n}}(x,y)$ is a Kummer extension of degree $m=(q^n+1)/(q+1)$.
			\item[$\bullet$] The full automorphism group ${{\rm Aut}}(\GK_{2,n})$ of $\GK_{2,n}$ acts on the set of $\mathbb{F}_{q^2}$-rational points of $\GK_{2,n}$ with two orbits, say $O_1$ and $O_2$, with \[O_1:= \{P_1,\dots,P_{q+1}\},\] lying over the $q+1$ points at infinity of $\cH_q$, and \[O_2:=\{R_1,\dots,R_{q^3-q}\},\] lying over the remaining $\mathbb{F}_{q^2}$-rational points.
		\end{itemize}
	\end{proposition}
	
	Before proceeding with the investigation of Weierstrass semigroups at points in $O_1$ and $O_2$, we describe the divisor associated to some algebraic functions on $\GK_{2,n}$ that we will use in the following.
	
	\begin{lemma}\cite[Lemma 4.2]{BM} \label{z}
		Let $  P_1:=(1:-1:0:0)\in O_1 $. Then
		\begin{itemize}
			\item[$1$.] $ (z) = \displaystyle \sum_{\substack{P \in \GK_{2,n}(\mathbb{F}_{q^2}) \\  P \notin  \{ P_1,\dots,P_{q+1}\}}}{\hspace{-15pt}P}-(q^2-q)\sum_{i=1}^{q+1}P_i $
			\item[$2$.] $(dz)=(q^{n+1}-q^n-q^2+2q-2) \sum_{i=1}^{q+1}P_i=\frac{2g(\GK_{2,n})-2}{(q+1)}\sum_{i=1}^{q+1}P_i.$
		\end{itemize}
	\end{lemma}
	
	\section{The Weierstrass semigroup at points in $O_1$} \label{resO1}
	
	In this section we investigate the structure of the Weierstrass semigroup $H(P)$ for $P\in O_1=\{P_1,\dots,P_{q+1}\}$. Note that \[P_i:=(1:a_i:0:0) \mbox{ and } a_i^{q+1}=1.\] These are all the points at infinity of $\GK_{2,n}$. Furthermore, in the function fields extension $\mathbb{F}_{q^{2n}}(x,y,z)/\mathbb{F}_{q^{2n}}(x,y)$ we have $P_i|\Pbar_i$, with $i=1,\dots,q+1$ where $\{\Pbar_i\}_{i}$ denotes the $q+1$ points at infinity of $\cH_q$. We can take \[P_1:=(1:-1:0:0)\] as a representative of points in $O_1$ to compute $H(P)$ for every $P \in O_1$. Indeed, it is known that the structure of Weierstrass semigroups is invariant under the action of automorphism groups see \cite[Lemma 3.5.2]{Sti} and $O_1$ is an ${\rm Aut}(\GK_{2,n})$-orbit. Thus, $\Pbar_1:=(1:-1:0)$.
	
	\begin{lemma} \label{mq+k}
		For all $k=0,\dots,\frac{m-1}{q^2-q}$ we have $mq + k(q^2-q) \in H(P_1)$. 
	\end{lemma}
	\begin{proof}
		Let $\varrho=x+y$. In \cite[Lemma 3.1]{BM} the function fields extension $\mathbb{F}_{q^{2n}}(x,y)/\mathbb{F}_{q^{2n}}(\varrho)$ is shown to be an Artin–Schreier extension of degree $q$ and in \cite[Proposition 2.1]{BM} it is pointed out that $\mathbb{F}_{q^{2n}}(x,y,z)/\mathbb{F}_{q^{2n}}(x,y)$ is a Kummer extension of degree $m$. This implies that 	
		\begin{equation}\label{rho}
		(\varrho)=mqP^{\infty}_1-m\sum_{i=2}^{q+1}P_i^{\infty},
		\end{equation} 
		see \cite[Equation (3.1)]{BM}.
		By Lemma~\ref{z} we can define for all $k=0,\dots,\frac{m-1}{q^2-q}$ the rational function $\vartheta_k=\frac{z^k}{\varrho}$. Then we have:
		\begin{align*}
		(\vartheta_k)&= k\sum_{\substack{P \in \GK_{2,n}(\mathbb{F}_{q^2}) \\  P \notin  \{ P_1,\dots,P_{q+1}\}}}{\hspace{-15pt}P}- (q^2-q)k\sum_{i=1}^{q+1}P_i - qmP_1+m\sum_{i=2}^{q+1}P_i \\
		&=k\sum_{\substack{P \in \GK_{2,n}(\mathbb{F}_{q^2}) \\  P \notin  \{ P_1,\dots,P_{q+1}\}}} \hspace{-15pt} P+(m-k(q^2-q))\sum_{i=2}^{q+1}P_i -(mq- k(q^2-q))P_1.
		\end{align*}
		Note that
		\begin{itemize}
			\item[$ \bullet $] $k\sum\limits_{\substack{P \in \GK_{2,n}(\mathbb{F}_{q^2}) \\  P \notin  \{ P_1,\dots,P_{q+1}\}}} \hspace{-15pt} P \geq 0$;
			\item[$ \bullet $] $(m-k(q^2-q))\sum_{i=2}^{q+1}P_i \geq 0$ as $ k\leq \frac{m-1}{q^2-q}$.
		\end{itemize}
		Thus,
		\[(\vartheta_k)=E-(mq- k(q^2-q))P_1,\]
		where $E \geq0$. This concludes the proof.
	\end{proof}
	
	\begin{obs}\label{fundamental}
		From Proposition~\textsl{\ref{BMprop}}, the curve $\GK_{2,n}$ is $ \mathbb{F}_{q^{2n}} $-maximal for odd $n \geq 3$. The fundamental equation \cite[Page xix (ii)]{HKT} guarantees that if $\mathcal{X}$ is an $\mathbb{F}_{q^2}$-maximal curve and $P \in \cX$ then \[qP+ \Phi^2(P) \equiv (q+1)P_0, \] where $P_0 \in \mathcal{X}(\mathbb{F}_{q^{2}})$ and $\Phi $ is the Frobenius homomorphism ($\alpha \mapsto \alpha^q$).
		Thus, if $P$ and $Q$ are $ \mathbb{F}_{q^{2n}} $-rational points of $\GK_{2,n}$ then \begin{equation}\label{fequation}
		(q^n+1)P \equiv (q^n+1)Q,
		\end{equation} as also $ \Phi^{2n}(P)=P $. In particular this implies that $q^n+1$ is a non-gap for every $P \in \GK_{2,n}(\mathbb{F}_{q^{2n}})$.
		For more details see also {\normalfont \cite[Proposition 10.6]{HKT}}.
	\end{obs}
	
	The following result allows us to construct explicitly a rational function which realize $q^n+1$ as a non-gap.
	
	\begin{lemma}\label{qn+1}
		Let $\alpha:=\big(\frac{x-1}{\varrho}\big)$. Then \[(\alpha)=(q^n+1)(Q-P_1), \] where $Q:=(1:0:1:0)$.
	\end{lemma}
	\begin{proof}
		Let $\Qbar$ be the affine point $\Qbar=(1:0:1)$. Then clearly $\Qbar$ is $\mathbb{F}_{q^2}$-rational. The tangent line at  $\mathcal{H}_q$ in $\Qbar$ is $t : X-1=0$. This line meets $\mathcal{H}_q$ only in $\Qbar$, so the intersection multiplicity is $q+1$. Thus the divisor of $t$ is \[(t)_{\mathbb{F}_{q^{2n}}(x,y)}=(q+1)\Qbar- \sum_{i=1}^{q+1}\Pbar_i.\]
		Using that $\mathbb{F}_{q^{2n}}(x,y,z)/\mathbb{F}_{q^{2n}}(x,y)$ is a Kummer extension, we obtain: \[(t)=m(q+1)Q- m\sum_{i=1}^{q+1}P_i=(q^n+1)Q- m\sum_{i=1}^{q+1}P_i.\]
		Considering now the quotient given by the rational functions $t$ and $\varrho=x+y$, we get
		\begin{equation*}
		\begin{split}
		\bigg(\frac{x-1}{\varrho}\bigg)&=(q^n+1)Q-mP_1-m\sum_{i=2}^{q+1}P_i + m\sum_{i=2}^{q+1}P_i -mqP_1 \\
		&=(q^n+1)Q-(q^n+1)P_1.
		\end{split}
		\end{equation*}
	\end{proof}
	
	Our aim is now to prove the following theorem, which is the main result of this section.
	
	\begin{theorem} \label{primo}
		Let $ H_1:=\langle mq+i(q^2-q),q^n+1 | i=0, \dots, s \rangle $. Then $H_1=H(P_1)$.
	\end{theorem}
	
	From Lemma~\ref{mq+k} and Lemma~\ref{qn+1}, $ H_1 $ is contained in $  H(P_1) $, so we only need to show the other inclusion to prove Theorem~\ref{primo}. 
	
	\subsection{Computing the genus of $H_1$}
	
	We recall that for a numerical semigroup $ G \subset \mathbb{N} $,  the \textit{genus} of $G$ is $g(G):=|\mathbb{N} \setminus G|$.
	In order to prove $H_1=H(P_1)$, we can equivalently show that the two semigroups $H_1$ and $H(P_1)$ have the same genus. Indeed, this means that there are not elements in $ H(P_1) \setminus H_1$ and hence $H_1=H(P_1)$. 
	We start by recalling the definition of an important class of numerical semigroups, called \textit{telescopic semigroups}, see \cite{Tel}.
	
	Let $(a_1, \dots, a_k)$ be a sequence of positive integers with greatest common divisor equal to $1$. Define \[d_i:=\gcd(a_1, \dots, a_i) \quad \mbox{ and } \quad A_i:=\Big\{\frac{a_1}{d_i},\dots, \frac{a_i}{d_i}\Big\},\]
	for $ i = 1, \dots,k $. Let $d_0:=0$ and $G_i$ be the semigroup generated by $A_i$. If $\frac{a_i}{d_i}\in G_{i-1}$ for all $i=2,\dots,k$ then the sequence $(a_1, \dots, a_k)$ is \textit{telescopic}. A numerical semigroup is called \textit{telescopic} if generated by a telescopic sequence.

	From \cite[Proposition 5.35]{Tel} the genus of a telescopic semigroup $S$ generated by a telescopic sequence $(a_1, \dots, a_k)$ is
	\begin{equation}\label{gs}
	g(S)=\frac{1}{2}\bigg(1+\sum_{i=1}^{k}\bigg(\frac{d_{i-1}}{d_i}-1\bigg)a_i\bigg).
	\end{equation}
	
	In order to compute the genus of $H_1$ we will proceed according to the following steps.
	
	\begin{itemize}
		
		\item[$ \bullet $] We construct a telescopic semigroup $S \subset H_1$;
		
		\item[$ \bullet $] we use Equation (\ref{gs}) to compute the genus of $S$, so that $g(H_1) \leq g(S) $; 
		
		\item[$ \bullet $] we compute explicitly the elements in $H_1 \setminus S$ and hence their number. In this way we get $g(H_1) \leq g(S)-|H_1 \setminus S| \leq g$. Since $g(H_1) \geq g$ as $H_1 \subseteq H(P_1)$ we get that $g(H_1)=g$. 
	\end{itemize}
	
	Let $S:=\langle mq, mq+(q^2-q),q^n+1 \rangle $.
	
	\begin{lemma} \label{S}
		The numerical semigroup $ S $ is telescopic. In particular, \[ g(S)=\frac{1}{2}(q(m^2-3m+2)+q^2(m-1)-q^n+q^{n+1}). \]
	\end{lemma}
	\begin{proof}
		The sequence $(mq, mq+(q^2-q),q^n+1)$ has GCD equal to $ 1 $. Furthermore we have $ q^n+1=m(q+1) $ and $ A_2= \{m,m+q-1\} $, so we can apply the equation (\ref{gs}). Thus,
		\begin{equation*}
		\begin{split}
		g(S) &= \frac{1}{2}(1+mq(-1)+(mq+q^2-q)(m-1)+(q-1)(q^n+1)) \\
		&= \frac{1}{2}(1-mq+m^2q+mq^2-mq-mq+q^2+q+q^{n+1}+q-q^n-1) \\
		&= \frac{1}{2}(q(m^2-3m+2)+q^2(m-1)-q^n+q^{n+1}).
		\end{split}
		\end{equation*}
	\end{proof}
	
	The following remark describes the elements in $H_1 \setminus S$.
	
	\begin{remark} \label{top}
		For every integer $ n $ there exist uniquely determined $ a $, $ b $ and $ c $ such that $ n=a(mq)+b(q^2-q)+c(q^n+1) $ , $ 0\leq b \leq m-1 $ and $ 0\leq c \leq q-1 $. In fact if $ a(mq)+b(q^2-q)+c(q^n+1)=a_1(mq)+b_1(q^2-q)+c_1(q^n+1) $ then we have $b \equiv b_1 $ mod $m$ and $c \equiv c_1 $ mod $q$.\\ 
		Moreover, $n \in S$ if and only if  $n=a(mq)+b(q^2-q)+c(q^n+1)$ with $ 0\leq b \leq m-1 $, $ 0\leq c \leq q-1 $ and $a \geq b$, because $ n=(a-b)(mq)+b(mq+q^2-q)+c(q^n+1) $.
	\end{remark}
	
	\begin{lemma}
		Let $n \in \mathbb{N}$ such that $n=a(mq)+b(q^2-q)+c(q^n+1) $, $ 0\leq b \leq m-1 $ and $ 0\leq c \leq q-1 $. Then $ n \in S $ if and only if $a \geq b$.
	\end{lemma}
	\begin{proof}
		If $a \geq b$, then $n=(a-b)(mq)+b(mq+q^2-q)+c(q^n+1) $, so $n \in S$.
		On the other hand, if $n \in S$ then $n=a(mq)+b(mq+q^2-q)+c(q^n+1)=(a+b)(mq)+b(q^2-q)+c(q^n+1) $.
	\end{proof}
	
	\begin{proposition}\label{contosa}
		For all $ i=1,\dots ,q^2-q-1 $ and $ j=q^2-q, \dots , (q^2-q)s-1 $ with $s=(m-1)/(q^2-q)$, define
		\[S_i=\lbrace i(mq)+(i+k_1)(q^2-q)+k_3(q^n+1) | k_1=1, \dots , is-i \mbox{ and }  k_3=0,\dots , q-1 \rbrace \]
		\[S_j=\{ j(mq)+(j+k_2)(q^2-q)+k_3(q^n+1) | 1\leq k_2\leq (q^2-q)s-j, \,  0\leq k_3 \leq q-1\}. \]
		Then
		\begin{itemize}
			\item[$1.$] $ S_i $, $ S_j $ $ \subset H_1 \setminus S $.
			\item[$2.$] $ \{S_i\}_i $ and $ \{S_j\}_j $ are families of mutually disjoint sets; we also have $ S_i \cap S_j = \emptyset$ for all $i$, $j$ in the corresponding ranges. 
			\item[$3.$] $|S_i|= (is - i)q$ and $|S_j| = ((q^2-q)s-j)q$. 
		\end{itemize}
	\end{proposition}
	
	\begin{proof} \vspace*{-5pt}\begin{itemize}
			\item[$1.$] If $ x \in S_i $ then $ x=i(mq)+(i+k_1)(q^2-q)+k_3(q^n+1)=i(mq+\frac{i+k_1}{i})(q^2-q)+k_3(q^n+1) $ with $\frac{i+k_1}{i} \leq s$, so $x \in H_1$. Moreover $  k_1 > 0 $ and we get $ x \notin S $. The same argument can be used for $S_j$.
			\item[$ 2. $] If $i(mq)+(i+k_1)(q^2-q)+k_3(q^n+1)=j(mq)+(j+k_2)(q^2-q)+k'_3(q^n+1)$, by Remark~\ref{top}, we have $i=j$, $k_3=k'_3$ and consequently $ k_1=k_2 $. In the same way it can be proved that the families are disjoint.
			\item[$ 3. $] From $1.$ and $2.$, we get the cardinality of $S_i$ (respectively $S_j$) simply multiplying the number of possibilities for $ k_1 $ (respectively $k_2$) by those for $ k_3 $.
		\end{itemize}
	\end{proof}	
	
	Picture \ref{fig1} describes the sets $S_i$ and $S_j$ for $q=2$ and $n=5$. In this picture a point of coordinates $(a,b)$ is used to represent the element $amq+b(q^2-q)+c(q^n+1)$ for some $0 \leq c \leq q-1$. Black dots represent elements of the numerical semigroup $S$, while white dots represent the elements contained in $S_i$ and $S_j$ for some $i$ and $j$.
	
	\begin{figure}[H]
		\label{fig1}
		\centering
		\begin{tikzpicture}[x=2.5cm,y=2.5cm]
		\draw[->] (0,0)--(3,0) node [right] {$a$};
		\draw[->] (0,0)--(0,2.9) node [above] {$b$};
		\draw[-] (0,2.6)--(3,2.6) node [right] {$10$};
		\draw[-] (0,0)--(2.7,2.7) node [above]{$a=b$};

		\foreach \y in {0}
		\foreach \x in {0,0.26,...,2.9}
		\node[draw,circle,inner sep=1pt,fill] at (\x,\y) {};
		\foreach \y in {0.26}
		\foreach \x in {0.26,0.52,...,2.9}
		\node[draw,circle,inner sep=1pt,fill] at (\x,\y) {};
		\foreach \y in {0.52}
		\foreach \x in {0.52,0.78,...,2.9}
		\node[draw,circle,inner sep=1pt,fill] at (\x,\y) {};
		\foreach \y in {0.78}
		\foreach \x in {0.78,1.04,...,2.9}
		\node[draw,circle,inner sep=1pt,fill] at (\x,\y) {};
		\foreach \y in {1.04}
		\foreach \x in {1.04,1.3,...,2.9}
		\node[draw,circle,inner sep=1pt,fill] at (\x,\y) {};
		\foreach \y in {1.3}
		\foreach \x in {1.3,1.56,...,2.9}
		\node[draw,circle,inner sep=1pt,fill] at (\x,\y) {};
		\foreach \y in {1.56}
		\foreach \x in {1.56,1.82,...,2.9}
		\node[draw,circle,inner sep=1pt,fill] at (\x,\y) {};
		\foreach \y in {1.82}
		\foreach \x in {1.82,2.08,...,2.9}
		\node[draw,circle,inner sep=1pt,fill] at (\x,\y) {};
		\foreach \y in {2.08}
		\foreach \x in {2.08,2.34,...,2.9}
		\node[draw,circle,inner sep=1pt,fill] at (\x,\y) {};
		\foreach \y in {2.34}
		\foreach \x in {2.34,2.6,...,2.9}
		\node[draw,circle,inner sep=1pt,fill] at (\x,\y) {};
		\foreach \y in {2.6}
		\foreach \x in {2.6,2.86,...,2.9}
		\node[draw,circle,inner sep=1pt,fill] at (\x,\y) {};
		
		
		
		\node[draw] at (3.2,1.04) {$S$};

		\foreach \x in {0.26}
		\foreach \y in {0.52,0.78,1.04}
		\node[draw,circle,inner sep=1pt] at (\x,\y) {};

		\foreach \x in {0.52}
		\foreach \y in {0.78,1.04,...,2.6}
		\node[draw,circle,inner sep=1pt] at (\x,\y) {};
		
		\foreach \x in {0.78}
		\foreach \y in {1.04,1.3,...,2.6}
		\node[draw,circle,inner sep=1pt] at (\x,\y) {};
		
		\foreach \x in {1.04}
		\foreach \y in {1.3,1.56,...,2.6}
		\node[draw,circle,inner sep=1pt] at (\x,\y) {};
		
		\foreach \x in {1.3}
		\foreach \y in {1.56,1.82,...,2.6}
		\node[draw,circle,inner sep=1pt] at (\x,\y) {};
		
		\foreach \x in {1.56}
		\foreach \y in {1.82,2.08,...,2.6}
		\node[draw,circle,inner sep=1pt] at (\x,\y) {};
		
		\foreach \x in {1.82}
		\foreach \y in {2.08,2.34,...,2.6}
		\node[draw,circle,inner sep=1pt] at (\x,\y) {};
		
		\foreach \x in {2.08}
		\foreach \y in {2.34,2.6}
		\node[draw,circle,inner sep=1pt] at (\x,\y) {};

		\foreach \x in {2.34}
		\foreach \y in {2.6}
		\node[draw,circle,inner sep=1pt] at (\x,\y) {};
		
		%
		%
		%
		%
		%
		
		%
		\draw[style=help lines,dashed] (0.26,0.52) --(0.26,2.6) node [above] { \textcolor{black}{$S_1$}};
		\draw[style=help lines,dashed] (0.52,0.78) --(0.52,2.6) node [above] { \textcolor{black}{$S_2$}};
		\draw[style=help lines,dashed] (0.78,1.04) --(0.78,2.6) node [above] { \textcolor{black}{$S_3$}};
		\draw[style=help lines,dashed] (1.04,1.3) --(1.04,2.6) node [above] { \textcolor{black}{$S_4$}};
		\draw[style=help lines,dashed] (1.3,1.56) --(1.3,2.6) node [above] { \textcolor{black}{$S_5$}};
		\draw[style=help lines,dashed] (1.56,1.82) --(1.56,2.6) node [above] { \textcolor{black}{$S_6$}};
		
		\draw[style=help lines,dashed] (1.82,2.08) --(1.82,2.6) node [above] {\textcolor{black}{$S_7$}};
		\draw[style=help lines,dashed] (2.08,2.34) --(2.08,2.6) node [above] {\textcolor{black}{$S_8$}};
		\draw[style=help lines,dashed] (2.34,2.6) --(2.34,2.6) node [above] { \textcolor{black}{$S_9$}};
		
		\end{tikzpicture}
		\caption{The sets $S_{i,j}$ and $S$ for $q=2$ and $n=5$}
		\label{figure}
		
		\begin{tabular}{r@{: }l r@{: }l}
			$\circ$ & Elements in $S_i$ and $S_j$ & \textbullet & Elements in $S$
		\end{tabular}
	\end{figure}
	
	We are in a position to prove our claim.
	
	\begin{theorem}
		We have $ g(H_1 ) = g $. In particular $H_1=H(P_1)$.
	\end{theorem}
	\begin{proof}
		By Proposition~\ref{contosa} we obtain the upper bound, \[g(H_1) \leq g(S) - \sum_{i=1}^{q^2-q-1} (is - i)q -\sum_{j=q^2-q}^{(q^2-q)s-1} ((q^2-q)s - j)q.\]
		Thus we get \begin{equation*} 
		\begin{split}
		g(H_1) & \leq g(S) - \sum_{i=1}^{q^2-q-1} (is - i)q -\sum_{j=q^2-q}^{(q^2-q)s-1} ((q^2-q)s - j)q  \\
		& = g(S) - \frac{1}{2} q (-1 - q + q^2) (-q + q^2) (-1 + s) - \sum_{k=1}^{(q^2-q)s-(q^2-q)} kq \\
		& =g(S) - \frac{1}{2} q (-1 - q + q^2) (-q + q^2) (-1 + s)+ \\
		& \quad -  \frac{1}{2} q (q - q^2 + (-q + q^2) s) (1 + q - q^2 + (-q + q^2) s) \\
		& =g(S) - \frac{1}{2} q (-1 - q + q^2) (-q + q^2) (-1 + s)+ \\
		& \quad -  \frac{1}{2} q (q^2-q)(s-1) (1 + (q^2-q)(s-1)) \\
		& = g(S) -\frac{1}{2} (-1 + q)^2 q^3 (-1 + s) s.
		\end{split}
		\end{equation*}
		Now, using that $s=\frac{m-1}{q^2-q}$ together with Lemma~\ref{S},
		\begin{equation*} 
		\begin{split}
		g(H_1) & \leq g(S) - \frac{1}{2}q(m-q^2+q-1)(m-1) \\
		& = \dfrac{1}{2} ((2 - 3 m + m^2) q + (-1 + m) q^2 - q^n + q^{n+1})\\ 
		&\quad- \frac{1}{2}q(m-q^2+q-1)(m-1).
		\end{split}
		\end{equation*}
		Finally,  \begin{equation*} \label{0}
		\begin{split} 
		g(H_1) &\leq \frac{1}{2}(mq^3-mq+q-q^3-q^n+q^{n+1}) \\
		&=\frac{1}{2}(\frac{q^{n+1}}{q+1}q(q^2-1)+q-q^3-q^n+q^{n+1}) \\
		&=\frac{1}{2}(q^{n+2}+q^2-q^3-q^n) \\
		&= \frac{1}{2}(-1 + q) (-q^2 + q^n + q^{1 + n}) = g. 
		\end{split}
		\end{equation*}
	\end{proof}
	
	\section{Weierstrass semigroup at points in $O_2$} \label{resO2}
	
	The second orbit of  $\mathbb{F}_{q^2}$-rational points of $\GK_{2,n}$ is $O_2=\GK_{2,n}(\mathbb{F}_{q^2}) / O_1=\{R_1,\dots,R_{q^3-q}\},$ that is, $O_2$ is the set of the $\mathbb{F}_{q^2}$-rational affine points of $\GK_{2,n}$. As already mentioned, the corresponding places are totally ramified in the function fields extension $\mathbb{F}_{q^{2n}}(x,y,z) / \mathbb{F}_{q^{2n}}(x,y)$; see Proposition \ref{ramificati}. 
	Let $\alpha$ be a given $(q+1) $-th root of -1 in $\mathbb{F}_{q^2}$ and let $R$ be the point \[R=(0:a:0:1).\] Note that  $R$ is in $O_2$. In fact $ a^{q^2-1}=a^{2(q+1)} = 1 $ and $R\in \GK_{2,n}$. Denote by $\Rbar$ the point 
	\[\Rbar=(0:a:1).\]
	Clearly, $\Rbar$ is an affine $\mathbb{F}_{q^2}$-rational point of the Hermitian curve $\cH_q$ and looking at the corresponding places, $R | \Rbar$ in the function fields extension $\mathbb{F}_{q^{2n}}(x,y,z)/\mathbb{F}_{q^{2n}}(x,y)$. As for the previous case, we can take $R$ to be a representative in $O_2$ because $H(R)=H(Q)$ for every $Q \in O_2$. 
	
	\begin{lemma}
		$q^n+1-k \in H(R)$ for all $k=0,\ldots,\frac{m-1}{q^2-q}$.
	\end{lemma}
	\begin{proof}
		Let $t_{R}$ be the tangent line in $\bar R$ at $\cH_q$. Then $t_{R}: Y-a=0$,	\[(y-a)_{\mathbb{F}_{q^{2n}}(x,y)}=(q+1) \Rbar - \sum_{i=1}^{q+1} \Pbar_i,\]
		
		and hence the divisor of $y-a$ in $\mathbb{F}_{q^{2n}}(x,y,z)$ is
		\[(y-a)=m(q+1) R -m \sum_{i=1}^{q+1} P_i=(q^n+1) R -m \sum_{i=1}^{q+1} P_i.\] 
		
		Define the algebraic function 	\[f_k=\frac{z^k}{y-a},\] where $z$ is defined as in {Lemma~\ref{z}} and $k=0,\ldots,(m-1)/(q^2-q)$. We get,
		\begin{equation*}
		\begin{split}
		(f_k)&=k\sum\limits_{\substack{Q \in \GK_{2,n}(\mathbb{F}_{q^{2n}}) \setminus H \\ Q \ne R}} Q+kR- k(q^2-q)\sum_{i=1}^{q+1} P_i-(q^n+1) R +m \sum_{i=1}^{q+1} P_i\\
		&=E-(q^n+1-k)R,
		\end{split}
		\end{equation*}
		where $E$ is an effective divisor. The claim now follows.
	\end{proof}
	
	\begin{lemma}
		$mq=q^n+1-m \in H(R)$.
	\end{lemma}
	
	\begin{proof}
		Since $T^{q+1}+1 \in \mathbb{F}_{q^2}[T]$ has $q+1$ distinct zeros in $\mathbb{F}_{q^2}$, the function $x$ has $q+1$ affine zeros in $\mathbb{F}_{q^{2n}}(x,y)$, namely the points $(0:\alpha:1)$ with $\alpha^{q+1}=-1$. Thus,
		\[(x)_{\mathbb{F}_{q^{2n}}(x,y)}=\Rbar + \sum_{\alpha^{q+1}=-1, \alpha \ne a} \Rbar_{(0:\alpha:0:1)} -  \sum_{i=1}^{q+1} \Pbar_i,\]
		and hence
		\[(x)=m R + m \sum_{\alpha^{q+1}=-1, \alpha \ne a} R_{(0:\alpha:0:1)} - m \sum_{i=1}^{q+1} P_i.\]
		Let $f=x/(y-a)$. Then, 
		\begin{equation*}
		\begin{split}
		(f)&=m R + m \sum_{\alpha^{q+1}=-1, \alpha \ne a} R_{(0:\alpha:0:1)} - m \sum_{i=1}^{q+1} R_i^\infty-(q^n+1) R +m \sum_{i=1}^{q+1} R_i^\infty\\
		&= m \sum_{\alpha^{q+1}=-1, \alpha \ne a} R_{(0:\alpha:0:1)}-(q^n+1-m)R,
		\end{split}
		\end{equation*}
		proving the statement.
	\end{proof}
	
	We are going to prove the following theorem, which is the main result of this section.
	
	\begin{theorem} \label{secondo}
		We have \[H(R)=\langle q^n+1-m,q^n+1-k \mid k=0,\ldots, (m-1)/(q^2-q)\rangle, \] for every $R \in O_2$.
	\end{theorem}
	
	In order to prove Theorem~\ref{secondo}, we proceed in a different way with respect to Section \ref{resO1}. Indeed this time we will determine the set of gaps $ G(R) $ in $R$ instead of $H(R)$. Doing so, we will then show that $H(R)=\mathbb{N}\setminus G(R)$ is exactly the semigroup $\langle q^n+1-m,q^n+1-k \mid k=0,\ldots, (m-1)/(q^2-q)\rangle$. 
	
	\subsection{Holomorphic differentials and gaps}
	
	We first recall some basic facts about holomorphic differentials on algebraic curves and function fields. 
	A differential $w$ on a function field $\mathbb{K}(x,y)$ is said to be \textit{holomorphic} if $(w) \geq 0$, see \cite{Sti} for details. For a divisor $D$ we set $\Omega(D):=\{w | w \mbox{ differential and } (w) \geq D\}$, so $\Omega(0)$ denotes the set of all holomorphic differentials. It is known that $\Omega(0)$ is a $\mathbb{K}$-vector space, with $\dim \Omega(0)=g$, where $g=g(\mathbb{K}(x,y))$, see \cite[Remark 1.5.12]{Sti}. 
	
	The main ingredient that we use to compute $G(R)$ is given by the following result.
	
	\begin{proposition} \label{Propbestng}{\rm \cite[Corollary 14.2.5]{VS}.}
		The integer $n \in \mathbb{N}$ is a gap at a place $P$ of $\mathbb{K}(x,y)$ if and only if there exists an holomorphic differential $\omega$ in $\mathbb{K}(x,y)$ such that \[\omega \in \Omega((n-1)P) \quad \mbox{ and } \quad \omega \notin \Omega(nP).\]
		In other words, $n$ is a gap at $P$ if and only if there exists an holomorphic differential $\omega$ such that \[v_P(\omega)=n-1.\]
	\end{proposition}
	
	According to this proposition, since $|G(R)|=g(\GK_{2,n})=\dim \Omega(0)$, our aim is to construct a basis for $\Omega(0)$ made by a class of holomorphic differentials having all distinct evaluations in $R$.
	
	First, we show that the function $z$ is a separating variable in the function field of the curve $\GK_{2,n}$ so that every differential of $\mathbb{F}_{q^{2n}}(x,y,z)$ can be written as $fdz$ for some $f$ in $\mathbb{F}_{q^{2n}}(x,y,z)$, see \cite[Proposition 4.1.8 (a)]{Sti}.
	
	By Lemma~\ref{z} we know that \[(dz)=(q^{n+1}-q^n-q^2+2q-2) \sum_{i=1}^{q+1}P_i=\frac{2g(\GK_{2,n})-2}{(q+1)}\sum_{i=1}^{q+1}P_i\] and $z$ is a separating variable from \cite[Proposition 4.1.8 (c)]{Sti}.
	
	We also observe that a differential $fdz$ is holomorphic if and only if for every $P \in \GK_{2,n}$ we have $v_P(f) \geq v_P(dz)$. Thus, by definition, $fdz$ is holomorphic if and only if \[f \in \mathscr{L}(dz)=\mathscr{L}((q^{n+1}-q^n-q^2+2q-2) \sum_{i=1}^{q+1}P_i),\] where $\mathscr{L}(dz)$ denotes the Riemann-Roch space of the divisor $(dz)$. Then we can construct an holomorphic differential $fdz$ provided that $f$ has only $P_i$ ($i=1\dots,q+1$) as poles, with multiplicity less than or equal to $q^{n+1}-q^n-q^2+2q-2$. 
	
	From Section \ref{resO1} we have 
	
	\begin{equation*}
	\begin{split}
	(z) &= \displaystyle \sum_{\substack{P \in \GK_{2,n}(\mathbb{F}_{q^2}) \\  P \notin  \{ P_1,\dots,P_{q+1}\}}}{\hspace{-15pt}P}-(q^2-q)\sum_{i=1}^{q+1}P_i ,\\
	(x)&=m R + m \sum_{\alpha^{q+1}=-1, \alpha \ne a} R_{(0:\alpha:0:1)} - m \sum_{i=1}^{q+1} P_i,\\
	(y-a)&=m(q+1) R -m \sum_{i=1}^{q+1} P_i=(q^n+1) R -m \sum_{i=1}^{q+1} P_i.
	\end{split}
	\end{equation*}
	
	Consider the family of differentials 
	
	\begin{equation*}
	\begin{split}
	\mathrm{I}:=\{&f_{j,k,l} dz  |\,  0 \leq k \leq m-1, \ 0 \leq l \leq q, \ 0 \leq j \leq q^2-2 \mbox{ with } \\ &k(q^2-q)+(j+l)m \leq q^{n+1}-q^n-q^2+2q-2\},
	\end{split}
	\end{equation*}
	where \[f_{j,k,l}=z^k (y-a)^j x^l.\]
	
	\begin{lemma} 
		The family $\mathrm{I}$ satisfies the following properties.
		\begin{itemize}
			\item[$ 1. $] $\mathrm{I} \subseteq \Omega(0)$.
			\item[$ 2. $] Let $f_{j,k,l}dz \in \mathrm{I}$. Then $v_R(f_{j,k,l}dz)=v_R(f_{j,k,l})=k+(q^n+1)j +l m$. 
			\item[$ 3. $] If $f_{j,k,l}dz,  f_{\tilde j, \tilde k, \tilde l} \in \mathrm{I}$ with $v_R(f_{i,j,k} dz)=v_R(f_{\tilde i, \tilde j, \tilde k} dz)$ then $(i,j,k)=(\tilde i, \tilde j, \tilde k)$.
		\end{itemize}
	\end{lemma}
	
	\begin{proof}
		\begin{enumerate}
			\item We want to show that  $f_{j,k,l} \in \mathscr{L}((q^{n+1}-q^n-q^2+2q-2) \sum_{i=1}^{q+1}R_i^\infty)$ for all $0 \leq k \leq m-1$, $0 \leq l \leq q$, $0 \leq j \leq q^2-2$ with $k(q^2-q)+(j+l)m \leq q^{n+1}-q^n-q^2+2q-2$.  We observe that with $j$, $k$ and $l$ fixed,
			\begin{equation*}
			\begin{split}
			(f_{j,k,l})&=(z^k (y-a)^j x^l)\\
			&=k\sum\limits_{\substack{Q \in \GK_{2,n}(\mathbb{F}_{q^{2n}}) \setminus H \\ Q \ne R}} Q+kR- k(q^2-q)\sum_{i=1}^{q+1} P_i+l m R+ \\
			&\:+ l m \! \! \! \sum\limits_{\substack{\alpha^{q+1}=-1 \\ \alpha \ne a}} R_{(0:\alpha:0:1)} - l m \sum_{i=1}^{q+1} R_i^\infty+j(q^n+1) R -jm \sum_{i=1}^{q+1} R_i^\infty\\
			&=E+(k+l m +j(q^n+1))R-(k(q^2-q)+l m +j m)\sum_{i=1}^{q+1} P_i,
			\end{split}
			\end{equation*}
			
			where $E$ is an effective divisor whose support is disjoint from $\{R,P_1,\ldots,P_{q+1}\}$.
			By definition, $f_{j,k,l} \in \mathscr{L}((q^{n+1}-q^n-q^2+2q-2) \sum_{i=1}^{q+1}P_i)$ if and only if $f_{j,k,l}$ has poles only in $\{P_i\}_i$ with multiplicity at most $q^{n+1}-q^n-q^2+2q-2$. 
			This is equivalent to require that $k(q^2-q)+l m +j m \leq q^{n+1}-q^n-q^2+2q-2$, which is satisfied by construction. We show that the bound considered for $j$, $k$, $l$ is not in contradiction with  the previous inequality. Actually we have \[\bigg\lfloor \frac{q^{n+1}-q^n-q^2+2q-2}{m}\bigg\rfloor=q^2-1-\bigg\lceil \frac{q^2-q+1}{m}\bigg\rceil=q^2-2,\]that is the greatest value for $j$ (note that it is strictly greater than the greatest value reached by $l$). Also, \[\bigg\lfloor \frac{q^{n+1}-q^n-q^2+2q-2}{q^2-q}\bigg\rfloor=q^{n-1}-1>m-1=q\frac{q^{n-1}-1}{q+1}.\]
			\item It follows directly by the computation of the divisor of $f_{i,j,k}$, as $v_R(dz)=0$.
			\item Suppose that
			\[k+(q^n+1)j +l m=v_R(f_{i,j,k} dz)=v_R(f_{\tilde i, \tilde j, \tilde k} dz)=\tilde k+(q^n+1)\tilde j +\tilde{l} m.\]
			By considering the above equality modulo $m$, and using that $k$ and $\tilde k$ are at most $m-1$ we obtain $k=\tilde k$. We need to show that \[(q+1)j +l=(q+1)\tilde j +\tilde{l}.\] As before, consider the equality modulo $q+1$. We have that $l$ and $\tilde{l}$ are less than or equal to $q$, so $l=\tilde{l}$. Clearly at this point $j=\tilde j$ and this proves the statement.
		\end{enumerate} 
	\end{proof}

	\begin{corollary}
		We have that \begin{equation*}
		\begin{split}
		L:=\{&k+(q^n+1)j +l m+1 | \,  0 \leq k \leq m-1, \ 0 \leq l \leq q, \ 0 \leq j \leq q^2-2 \mbox{ and}\\ &k(q^2-q)+(j+l)m \leq q^{n+1}-q^n-q^2+2q-2\} \subseteq G(R).
		\end{split}
		\end{equation*}
		
		Moreover, the differentials $f_{j,k,l} dz$ are linearly independent over $\mathbb{F}_{q^{2n}}$, as they have all distinct evaluations in $R$.
	\end{corollary}
	
	The following observation explains our next step.

	\begin{obs}\label{genL}
		If we show that for a subset $L \subseteq G(R)$, \[|L|=g(\GK_{2,n})=(q-1)(q^{n+1}+q^n-q^2)/2\] then	\begin{itemize}
			\item[$\bullet$] $L=G(R)$,
			\item[$\bullet$] $\mathrm{I}$ is a $\mathbb{F}_{q^{2n}}$-basis for the vector space $\Omega(0)$ of $\GK_{2,n}$.
		\end{itemize}
	\end{obs}
	
	According to Observation~\ref{genL}, we only need to show that $ L $ contains exactly $g(\GK_{2,n})$ elements. Doing so, the following theorem is proven.
	
	\begin{theorem} \label{main}
		Let $q$ be a prime power and $n \geq 5$ odd. Let also $R$ be an affine $\mathbb{F}_{q^2}$-rational point of the curve $\GK_{2,n}$. Let $L$ be the set
		\begin{equation*}
		\begin{split}
		L=\{&k+(q^n+1)j +l m+1 | \,  0 \leq k \leq m-1, \ 0 \leq l \leq q, \ 0 \leq j \leq q^2-2 \mbox{ with }\\ &k(q^2-q)+(j+l)m \leq q^{n+1}-q^n-q^2+2q-2\}.
		\end{split}
		\end{equation*}
		Then we have
		\[|L|=g(\GK_{2,n})=\frac{1}{2}(q-1)(q^{n+1}+q^n-q^2).\] In particular, $L=G(R)$.
	\end{theorem}
	
	In order to prove Theorem~\ref{main} we are going to proceed with a direct computation of the number of distinct elements in $ L $. First, we make clear the range of the entries of $(k,l,j)$ for elements in $L$.
	
	\begin{obs}
		Note that if $k$ could always reach the value $m-1$, the set $L$ would be the whole $\mathbb{N}$. As a matter of fact, we can move from $ h \in L $, with associated value $ (k,l,j) $, to $h+1$ simply putting $k=0$ and by replacing $l$ (or $ j $) with  $l+1$ (or $j+1$). This clearly cannot be possible by the bound $k(q^2-q)+(j+l)m \leq q^{n+1}-q^n-q^2+2q-2$. However when $k_{max}(j,l)$ coincides with $m-1$ then from $(q^n+1)j +l m $ to $ (q^n+1)j +(l+1) m$ there are not non-gaps.
	\end{obs}
	
	For the sake of simplicity, in the following we identify an element in $L$ with the associated triple $(k,l,j)$. Furthermore $k_{max}$ always denotes $k_{max}:=k_{max}(l+j)$.
	
	\begin{lemma}
		Let $(k,l,j)$ be in $ L$, then $j+l <q^2-1$ for all $k \geq 0$.
	\end{lemma}
	
	\begin{proof}
		Suppose that $j+l\geq q^2-1$. Then
		
		\begin{equation*}
		\begin{split}
		(j+l)m &\geq (q^2-1)m \\
		&= (q+1)(q-1)\frac{q^n+1}{q+1}\\
		&=(q-1) (q^n+1)\\
		&=q^{n+1}-q^n+q-1>q^{n+1}-q^n-q^2+2q-2.
		\end{split}
		\end{equation*} 
		In the last equation we used that $q^2-q+1>0$ for every integer $q$.
	\end{proof}
	
	\begin{lemma}\label{max1}
		Let $ (k,l,j) $ be in $  L $ with $j+l < q-1$. Then $k_{max}=m-1$.
	\end{lemma}
	\begin{proof}
		By direct computation,
		\begin{equation*}
		\begin{split}
		(m-1)(q^2-q)+(q-2)m &= mq^2-q^2-mq+q+mq-2m\\
		&=m(q^2-2)-q(q-1)\\
		&<m(q^2-1)-q(q-1)=(q^n+1-q)(q-1)\\
		&=q^{n+1}-q^n-q^2+2q-1.
		\end{split}
		\end{equation*}
		
		Hence $(m-1)(q^2-q)+(q-2)m \leq q^{n+1}-q^n-q^2+2q-2$ and the claim follows.
	\end{proof}
	
	\begin{lemma} \label{max2}
		Let $(k,l,j)$ be in $L$. Then $k_{max}=m-1-(s(j+l-q+1)+1)$ for all $q-1 \leq j+l < q^2-1$.
	\end{lemma}
	\begin{proof}
		Let $\bar{k}=m-1-(s(j+l-q+1)+1) $. We note that if $j+l \geq q-1$ then $\bar{k} \leq m-2$. If on the contrary $j+l < q^2-1$, then
		\begin{equation*}
		\begin{split}
		\bar{k}(q^2-q)+(j+l)m &< (m-1-(s(q^2-1-q+1)+1))(q^2-q)\!+\!(q^2-1)m\\
		&= (m-2-s(q^2-q))(q^2-q)+(q^2-1)m\\
		&=(m-2-m+1)(q^2-q)+m(q+1)(q-1)\\
		&=-q(q-1)+(q^n+1)(q-1)\\
		&=(q-1)(q^n+1-q)=q^{n+1}-q^n-q^2+2q-1.
		\end{split}
		\end{equation*}
		
		So we prove that $\bar{k}(q^2-q)+(j+l)m \leq q^{n+1}-q^n-q^2+2q-2$. Using that $j+l \geq q-1$, we get \begin{equation*}
		\begin{split}
		(\bar{k}+1)(q^2-q)+(j+l)m &\geq (m-1-s(q-1-q+1))(q^2-q)+(q-1)m\\
		&= (m-1)(q^2-q)+m(q-1)\\
		&=m(q^2-1)-q(q-1)\\
		&=(q^n-q+1)(q-1)=q^{n+1}-q^n-q^2+2q-1,
		\end{split}
		\end{equation*}
		
		and thus $(\bar{k}+1)(q^2-q)+(j+l)m >  q^{n+1}-q^n-q^2+2q-2$, proving that $\bar{k}=k_{max}$.
	\end{proof}
	
	\begin{obs}
		The value $k_{max}$ depends on the sum of $j$ and $l$, i.e. $k_{max}(j,l)=k_{max}(j+l)$.
	\end{obs}
	
	We are now in a position to prove Theorem~\ref{main}.
	
	\noindent\begin{proof}[Proof, Theorem~\ref{main}]
		
		Let $M_q=q^{n+1}-q^n-q^2+2q-2$. Then $L$ can be described as follows 
		
		\begin{equation*}\begin{split}
		L:=\{&k+(q^n+1)j +l m+1 |\,   0 \leq k \leq k_{max}, \ 0 \leq l \leq q, \ 0 \leq j+l \leq q^2-2 \mbox{ :}\\ 
		&k(q^2-q)+(j+l)m \leq M_q\}.
		\end{split}
		\end{equation*}
		
		From Lemma~\ref{max1} and Lemma~\ref{max2} we can partition $L$ into two subsets
		
		$$L_1:=\{k+(q^n+1)j +l m+1 \mid  0 \leq k \leq m-1, \ 0 \leq j+l \leq q-2 \ : \ k(q^2-q)+(j+l)m \leq M_q\}$$
		and 
		$$L_2:=\{k+(q^n+1)j +l m+1 \mid  0 \leq k \leq k_{max}, \ 0 \leq l \leq q,$$
		$$ \ q-1 \leq j+l \leq q^2-2 \ : \  k(q^2-q)+(j+l)m \leq M_q\}.$$
		
		Since every element of $L$ is uniquely determined by $(k,l,j)$, we just need to compute the number of triples which are contained in $L_1$ and $L_2$. \\
		
		\begin{itemize}
			\item[$\bullet$] \textit{The computation of $|L_1|$.} We note that $k_{max}$ in $L_1$ does not depend on $l$ and $j$. Let $l$ be fixed. Then we have $m$ choices for $k$, while $j=0,\dots,q-2-l$. Thus,
			
			\[ |L_1|=\sum_{l=0}^{q-2}{m(q-l-1)}=m\frac{q(q-1)}{2}. \]
			
			\item[$\bullet$] \textit{The computation of $|L_2|$.} We distinguish two cases, namely $l\leq q-1$ and $l=q$. If $l\leq q-1$, then $j=q-1-l,\dots,q^2-2-l$; while in the latter case $j=0,\dots,q^2-2-l$. Suppose now that $(l,j)$ is fixed. Repeating the previous argument we note that at each step $k=0,\dots,m-2-s(j+l-q+1)$; so we have $  m-1-s(j+l-q+1) $ choices for $k$ . By adding up for every value of $l$ and $j$, we obtain
			\begin{equation*}
			\begin{split}
			|L_2|&=\sum_{l=0}^{q-1}\sum_{j+l=q-1}^{q^2-2}\!{\!(m-1-s(j+l-q+1))}\!+\!\!\!\!\sum_{j=0}^{q^2-q-2}{(m-1-s(j+1))}\\
			&=\sum_{l=0}^{q-1}((q^2-q)(m-1)-s\frac{(q^2-q)(q^2-q-1)}{2})+\\
			&\quad+\sum_{j=0}^{q^2-q-2}{(m-1-s(j+1))}\\
			&=\frac{1}{2}q(q^2-q)(2(m-1)-s(q^2-q-1))\!+\!\!\sum_{j=0}^{q^2-q-2}\!{\!(m-1-s(j+1))}\\
			&=\frac{1}{2}q(q^2-q)(2(m-1)-s(q^2-q-1))+((q^2-q-1)(m-1)+\\
			&\quad-\frac{1}{2}s(q^2-q)(q^2-q-1))\\
			&=\frac{1}{2}q(q^2-q)(2(m-1)-s(q^2-q-1))+\\
			&\quad+\frac{1}{2}(q^2-q-1)(2(m-1)-s(q^2-q)).
			\end{split}
			\end{equation*}
			
			Putting $ m-1=s(q^2-q) $ in the equation we have
			\begin{equation*}
			\begin{split}
			|L_2|&=\frac{1}{2}q(q^2-q)s(2q^2-2q-q^2+q+1)+\frac{1}{2}(q^2-q-1)s(q^2-q)\\
			&=\frac{1}{2}q(q^2-q)s(q^2-q+1)+\frac{1}{2}(q^2-q-1)s(q^2-q)\\
			&=\frac{1}{2}s(q^2-q)(q^3-q^2+q+q^2-q-1)=\frac{1}{2}(m-1)(q^3-1).
			\end{split}
			\end{equation*}
		\end{itemize}
		Hence,
		\begin{equation*}
		\begin{split}
		|L|=|L_1|+|L_2|&=\frac{1}{2}mq(q-1)+\frac{1}{2}(m-1)(q-1)(q^2+q+1)\\
		&=\frac{1}{2}(q-1)(mq+mq^2+mq+m-q^2-q-1)\\
		&=\frac{1}{2}(q-1)(m(q+1)+m(q)(q+1)-q^2-q-1)\\
		&=\frac{1}{2}(q-1)(q^n+1+(q^n+1)(q)-q^2-q-1)\\
		&=\frac{1}{2}(q-1)(q^{n+1}+q^n-q^2)=g(\GK_{2,n}),
		\end{split}
		\end{equation*}
		proving the statement.
	\end{proof}
	
	\section{On the Frobenius dimension of $\GK_{2,n}$} \label{secFrob}
	
	In this section we investigate the Frobenius dimension of the curves $\GK_{2,n}$. In particular, we are interested in comparing it with the Frobenius dimension of first generalized GK curve $\GK_{2,n}$.
	
	Recall that for any $\mathbb{F}_{q^{2n}}$-maximal curve $\cX_n$, the Fundamental Equation \eqref{fequation} is written as \[q^nP+\Phi^{2n}(P)\sim(q^n+1)P_0,\] where $P\in \cX_n$, $P_0 \in \cX_n(\mathbb{F}_{q^{2n}})$ and $\Phi$ is the Frobenius automorphism. 
	
	The complete linear series given by $\cD:=|(q^n+1)P_0|$ is said to be the \textit{Frobenius linear series} of $\cX_n$ and $r:=\dim \cD$ is the \textit{Frobenius dimension} of $\cX_n$. This dimension is one of the most important (birational)  invariants of a maximal curve. 
	
	The following proposition allows us to easily compute the Frobenius dimension of any maximal curve, see \cite[Section 10.2]{HKT}.
	
	\begin{proposition}{\cite[Propositions 10.6 and 10.9]{HKT}} \label{frob}
		Let $\cX_n$ be an $\mathbb{F}_{q^{2n}}$-maximal curve having Frobenius dimension $r$. Let $P$ be an $\mathbb{F}_{q^{2n}}$-rational point of $\cX_n$. Then the following holds.
		\begin{itemize}
			\item[$1.$] The Frobenius dimension $r$ coincides with the number of non-trivial non-gaps at $P$ which are less than or equal to $q^n$, i.e. \[0<m_1(P)<\cdots<m_{r-1}(P)\leq q<m_r(P). \]
			\item[$2.$] $r\geq 2$ holds. If $r \geq 3$ and $m_{r-2}(P)<q^n-1$, then $P$ is a Weierstrass point of $\cX_n$.
		\end{itemize}
	\end{proposition}
	
	The following theorem is obtained using the results of Sections \ref{resO1} and \ref{resO2} together with Proposition \ref{frob}.
	
	\begin{theorem}\label{frobx}
		Let $q$ be a prime power and $n\geq 5$ an odd integer. The Frobenius dimension of the second generalized $GK$ curves $\GK_{2,n}$ is
		
		\begin{equation*}\label{rdim}
		r=\frac{m-1}{q^2-q}+2, \mbox{ where } m=\frac{q^n+1}{q+1}.
		\end{equation*}
	\end{theorem}
	
	We now show two applications of Theorem \ref{rdim}. 
	
	It is known that if $H(P)$ is symmetric for $P \in \GK_{2,n}$, that is, $2g-1 \in G(P)$, then $P$ is a Weierstrass point of $\GK_{2,n}$; see \cite[Proposition 50]{Konto}. However, the converse is not necessarily true. The results obtained in Sections \ref{resO1} and \ref{resO2} together with Theorem \ref{rdim} allow us to provide a counterexample.
	
	\begin{corollary} \label{wp}
		If $P \in \GK_{2,n}(\mathbb{F}_{q^{2n}})$ then $P$ is a Weierstrass point of $\GK_{2,n}$.
		In particular, if $P \in O_1 \cup O_2$ then $H(P)$ is not symmetric even though $P$ is a Weierstrass point of $\GK_{2,n}$.
	\end{corollary}
	
	\begin{proof}
		The fact that $H(P)$ is not symmetric follows from a direct computation from Theorems \ref{primo} and \ref{secondo}. The claim regarding the points in $\GK_{2,n}(\mathbb{F}_{q^{2n}})$ follows from \cite[Corollary 4.5.]{FT} and Theorem \ref{rdim} noting that $q^n+1-r<g(\GK_{2,n})$.
	\end{proof}

	Theorem \ref{frobx} has also the following second application. It allows us to exhibit another way to prove that the curves $\GK_{2,n}$ and $\GK_{1,n}$ are not isomorphic for any $n \geq 5$ odd. Indeed in \cite[Corollary 2.6]{BM}, the authors proved the following theorem.
	
	\begin{theorem}\label{isomor}
		Let $q$ be a prime power and $n \geq 3$ be odd. Then $ \GK_{2,n}$ is isomorphic to $\GK_{1,n}$ if and only if $n =3$.
	\end{theorem} 
	
	In order to prove Theorem~\ref{isomor} Beelen and Montanucci made use of the theory of automorphism groups of algebraic curves. In fact, since the full automorphism group of an algebraic curve is invariant under (birational) isomorphisms, it is sufficient to note that $ {\rm Aut}(\GK_{1,n}) \ne {\rm Aut}(\GK_{2,n})$, see \cite{BM} and \cite{GMP, GOS}. 
	
	We instead are going to compute the Frobenius dimension of the curves $\GK_{1,n}$ and compare it with the value obtained in Theorem~\ref{frobx}.
	
	Remember that the function field of $\GK_{1,n}$ is the compositum of the function fields $ \mathbb{F}_{q^{2n}}(x,y) $ and $\mathbb{F}_{q^{2n}}(y,z)$, where $x,y$ and $z$ satisfy Equation \eqref{GGS}. Let $P_\infty$ denote the common pole of $x,y$ and $z$. In \cite{GOS}, the following result is proved.
	
	\begin{proposition}
		The set of non-gaps at $P_\infty$ in $\GK_{1,n}$ is \[ \{i(q^n+1)+jmq+kq^3 | i,j,k \in \mathbb{N} \mbox{ with }0\leq i \leq q, \, 0 \leq j \leq q^2, \, k\geq 0 \}. \]
	\end{proposition}
	
	As for the curve $\GK_{2,n}$ we can now compute the Frobenius dimension of $\GK_{1,n}$.
	
	\begin{corollary} \cite[Corollary 3.43]{Pietro}
		The Frobenius dimension of the $\GK_{1,n}$ curve is equal to 
		\begin{equation*}
		r^\prime = q^{n-3}+ \sum_{i=2}^{n-2}(-1)^{i+1}q^i+1.
		\end{equation*}
	\end{corollary}
	
	Finally, we show how this computation allows us to obtain another proof of Theorem \ref{isomor}.
	
	\begin{theorem}
		Let $q$ denote a prime power and $n \geq 3$ an odd integer. For a fixed $n$, let $\GK_ {1,n}$ and $\GK_{2,n}$ be the first and the second generalized $GK$ curve respectively. Then $\GK_ {1,n}$ and $\GK_ {2,n}$ are not isomorphic for every $n \geq 5$.
	\end{theorem}
	
	\begin{proof}
		The proof follows directly by comparing the Frobenius dimensions of $\GK_ {1,n}$ and $\GK_ {2,n}$. Indeed, we have \begin{equation*}
		\begin{split}
		q^{n-3}+ \sum_{i=2}^{n-2}(-1)^{i+1}q^i+1 \geq q^{n-3}+q^{n-2}-q^{n-3}+1,
		\end{split}
		\end{equation*}
		
		and we want to show that
		
		\begin{equation}\label{conti}
		\begin{split}
		q^{n-3}+q^{n-2}-q^{n-3}+1 &\geq \frac{m-1}{q^2-q}+2\\
		&=\dfrac{\frac{q^n+1}{q+1}-1}{q^2-q}+2\\
		&=\frac{q^{n-1}-1}{q^2-1}+2.
		\end{split}
		\end{equation}
		
		However, we obtain equation \eqref{conti} by computing
		\begin{equation}
		\begin{split}
		&q^{n-2}+1- \frac{q^{n-1}-1}{q^2-1}-2=\\
		=& (q^{n-2}-1)(q^2-1)-q^{n-1}+1=\\
		=& q^{n}-q^{n-2}-q^{n-1}-q^2+2,
		\end{split}
		\end{equation}
		which is a positive integer for $n \geq 5$.
	\end{proof}
	
	\section{Applications to AG codes} \label{appAG}
	
	In this section we will apply the results obtained in Sections \ref{resO1} and \ref{resO2} to construct AG codes and AG quantum codes from the second generalized GK curve. Explicit tables containing the parameters of the resulting codes for  $q=2$ and $n=5$ can be found in Subsection \ref{tables}.
	
	As before, $\GK_{2,n}(\mathbb{F}_{q^{2n}})$ denotes the set of all $\mathbb{F}_{q^{2n}}$-rational points of $\GK_{2,n}$ while $\mathbb{F}_{q^{2n}}(\GK_{2,n})$ denotes the set of $\mathbb{F}_{q^{2n}}$-rational functions on $\GK_{2,n}$. A divisor $D$ is $\mathbb{F}_{q^{2n}}$-rational if it is fixed by the Frobenius endomorphism $\Phi^{2n}$.
	
	We briefly recall the definition of an AG code, see \cite[Chapter 2]{Sti} and \cite{Tel} for a more detailed description. 
	Let $P_1, \dots P_N \in \GK_{2,n}(\mathbb{F}_{q^{2n}})$ be pairwise distinct points and consider the divisor $D:=P_1+\cdots P_N$. Let $G$ be another $\mathbb{F}_{q^{2n}}$-rational divisor whose support is disjoint from the one of $D$. Let $e$ denote the following linear map
	
	$$\begin{cases} e \colon \mathscr{L}(G) \rightarrow \mathbb{F}_{q^{2n}}^N, \\ \alpha \mapsto e(\alpha):= (\alpha(P_1),\dots,\alpha(P_N)). \end{cases}$$
	
	The AG code associated to $D$ and $ G $ is $C(D,G):=e(\mathscr{L}(G))$. The code $C(D,G)$ is an $ [N,k,d]_{q^{2n}}$-code with $d \geq N-\deg G$ and $k=\ell(G)-\ell(G-D)$. When $G:=n_P P$, $n_P \in \mathbb{N}$, then $C(D,G)$ is a \textit{one-point code}.
	The dual code $C^{\perp}(D,G)$ is an AG code with dimension $k^\perp:=N-k$ and minimum distance $d^\perp \geq \deg G -2g+2$, where $g$ is the genus of $\GK_{2,n}$.
	
	Let $P \in \GK_{2,n}(\mathbb{F}_{q^{2n}})$ and set \[H(P)=\{\varrho_1:=0<\varrho_2<\cdots \}\] the Weierstrass semigroup at $P$. For $\ell >0$ the \textit{Feng-Rao function} is defined as \[\nu_\ell :=|\{(i,j) \in \mathbb{N}^2 : \varrho_i+\varrho_j=\varrho_\ell\} |.\] Consider now the AG code $C_\ell:=C^\perp(P_1+\cdots+P_N, \varrho_\ell P)$, with $N >\ell+1$ and $P_1,\dots,P_N, P$ distinct points.
	
	\begin{proposition} \cite[Theorem 5.24]{Tel}\label{feng}
		$C_\ell$ is a linear $[N,k,d]_{q^{2n}}$-code with $k=N-\ell$ and $d\geq d_{ORD}(C_\ell)$, where \[d_{ORD}(C_\ell):=\min \{\nu_m | m\geq \ell\}\] is the Feng-Rao designed minimum distance.
	\end{proposition}
	
	From \cite[Theorem 5.24]{Tel} we have also the following proposition which shows that for large values in $H(P)$, $d_{ORD}$ can be easily computed.
	
	\begin{proposition} \label{limitcase}
		Let $ H(P) $ be a Weierstrass semigroup. Then $ d_{ORD}(C_\ell(P)) \geq \ell+1-g $ and
		equality holds if $\varrho_{\ell}+1 \geq 4g$.
	\end{proposition}
	
	\subsection{Tables of AG codes} \label{tables}
	
	Here we are going to show the tables of AG dual codes constructed on the curve $\GK_{2,n}$. We consider just the case $q:=2$ and $n:=5$. Consider the set $\GK_{2,5} (\mathbb{F}_{2^{10}})$, with $|\GK_{2,5}(\mathbb{F}_{2^{10}})|=3969$. If we take a point $P \in O_1$, then we have \[H(P):=\langle mq+i(q^2-q),q^n+1 | i=0, \dots, s \rangle, \] by Theorem \ref{primo}. Thus, we can calculate the parameters of $C^1_\ell:=C^\perp(P_1+\cdots+P_{3968}, \varrho_\ell P)$. The following table is made using the software Magma \cite{MAGMA}.

	\begin{table}[H]
		
		\centering
		\begin{normalsize}

			\begin{tabular}{|c||c|c|c|c||c|c|c|c||c|c|c|c|}
				\hline
				$ n $ & $ k $ & $ \varrho_{\ell} $ & $ \nu_\ell $ &  $ d_{ORD} $ & $ k $ & $ \varrho_{\ell} $ & $ \nu_\ell $ &  $ d_{ORD} $ & $ k $ & $ \varrho_{\ell} $ & $ \nu_\ell $ &  $ d_{ORD} $\\
				\hline \hline
				3968&3967&0&1&1& 3966&22&2&2& 3965&24&2&2\\ \hline
				3968&3964&26&2&2& 3963&28&2&2& 3962&30&2&2\\ \hline
				3968&3962&30&2&2& 3961&32&2&2& 3960&33&2&2\\ \hline
				3968&3959&44&3&3& 3958&46&4&3& 3957&48&5&3\\ \hline
				3968&3956&50&6&3& 3955&52&7&3& 3954&54&8&3\\ \hline
				3968&3953&55&4&3& 3952&56&7&3& 3951&57&4&3\\ \hline
				3968&3950&58&6&3& 3948&60&5&3& 3947&61&4&3\\ \hline
				3968&3946&62&4&3& 3945&63&4&3& 3944&64&3&3\\ \hline
				3968&3943&65&4&4& 3942&66&5&5& 3941&68&6&6\\ \hline
				3968&3940&70&8&6& 3939&72&10&6& 3938&74&12&6\\ \hline
				3968&3937&76&14&6& 3936&77&6&6& 3935&78&14&6\\ \hline
				3968&3934&79&8&6& 3933&80&14&6& 3932&81&10&6\\   \hline 
				39688&3931&82&14&6& 3930&83&12&6&3929&84&14&6\\   \hline 
			\end{tabular} 
			
		\end{normalsize}	
		
		\caption{Parameters of $C^1_\ell$.}
		\label{tab2}
	\end{table}
	
	\begin{table}
		\centering
		\begin{normalsize}

			\begin{tabular}{|c||c|c|c|c||c|c|c|c||c|c|c|c|}
				\hline
				$ n $ & $ k $ & $ \varrho_{\ell} $ & $ \nu_\ell $ &  $ d_{ORD} $ & $ k $ & $ \varrho_{\ell} $ & $ \nu_\ell $ &  $ d_{ORD} $ & $ k $ & $ \varrho_{\ell} $ & $ \nu_\ell $ &  $ d_{ORD} $\\
				\hline \hline
				39688&3928&85&14&6&3927&86&14&6&3926&87&16&6\\   \hline 
				39688&3925&88&17&6&3924&89&14&6&3923&90&18&6\\   \hline 
				39688&3922&91&12&6&3921&92&19&6&3920&93&10&6\\   \hline 
				39688&3919&94&20&6&3918&95&8&6&3917&96&21&6\\   \hline 
				39688&3916&97&6&6& 3915&98&22&8&3914&99&8&8\\   \hline 
				39688&3913&100&21&12&3912&101&12&12&3911&102&22&16\\   \hline 
				39688&3910&103&16&16&3909&104&23&20&3908&105&20&20\\   \hline 
				39688&3907&106&24&24&3906&107&24&24&3905&108&25&25\\   \hline 
				39688&3904&109&28&28&3903&110&29&28&3902&111&28&28\\   \hline 
				39688&3901&112&31&28&3900&113&28&28&3899&114&33&28\\   \hline 
				39688&3898&115&28&28&3897&116&35&28&	3896&117&28&28\\   \hline 
				3968&3895&118&37&28&3894&119&28&28&3893&120&39&30\\ \hline
				3968&3892&121&30&30&3891&122&39&32&3890&123&32&32\\ \hline
				3968&3889&124&39&34&3888&125&34&34&3887&126&39&36\\ \hline
				3968&3886&127&36&36&3885&128&39&38&3884&129&38&38\\ \hline
				3968&3883&130&39&39&3882&131&40&40&3881&132&41&41\\ \hline
				3968&3880&133&42&42&3879&134&44&44&3878&135&44&44\\ \hline
				3968&3877&136&47&46&3876&137&46&46&3875&138&50&48\\ \hline
				3968&3874&139&48&48&3873&140&53&50&3872&141&50&50\\ \hline
				3968&3871&142&56&52&3870&143&52&52&3869&144&57&54\\ \hline
				3968&3868&145&54&54&3867&146&58&56&3866&147&56&56\\ \hline
				3968&3865&148&59&58&3864&149&58&58&3863&150&60&60\\   \hline 
				3868&3862&151&60&60&3861&152&61&61&3860&153&62&62\\   \hline 		
				3868&3859&154&63&63&3858&155&64&64&3857&156&65&65\\   \hline 
				3868&3856&157&66&66&3855&158&67&67&3854&159&68&68\\   \hline 
				3868&3853&160&69&69&3852&161&70&70&3851&162&71&71\\   \hline 
				3868&3850&163&72&72&3849&164&73&73&3848&165&74&74\\   \hline 
				3868&3847&166&75&75&3846&167&76&76&3845&168&77&77\\   \hline 
				3868&3844&169&78&78&3843&170&79&79&3842&171&80&80\\   \hline 
				3868&3841&172&81&81&3840&173&82&82&3839&174&83&83\\   \hline 
				3868&3838&175&84&84&3837&176&85&85&3836&177&86&86\\   \hline 
				3868&3835&178&87&87&3834&179&88&88&3833&180&89&89\\   \hline 
				3868&3832&181&90&90&3831&182&91&91&3830&183&92&92\\   \hline 
			\end{tabular}
			
		\end{normalsize}	
		
		\caption{Parameters of $C^1_\ell$.}
		
	\end{table} 
	
	Consider now a point $R \in O_2$. From Theorem \ref{secondo} we know that \[H(R)=\langle q^n+1-m,q^n+1-k \mid k=0,\ldots, (m-1)/(q^2-q)\rangle.\] Applying the same argument as above, we can compute  the parameters of the AG code $C^2_\ell:=C^\perp(P_1+\cdots+P_{3968}, \varrho_\ell R)$, when $q:=2$ and $n:=5$.
	
	\begin{table}
		\centering
		
		\begin{tabular}{|c||c|c|c|c||c|c|c|c||c|c|c|c|}
			\hline
			$ n $ & $ k $ & $ \varrho_{\ell} $ & $ \nu_\ell $ &  $ d_{ORD} $ & $ k $ & $ \varrho_{\ell} $ & $ \nu_\ell $ &  $ d_{ORD} $ & $ k $ & $ \varrho_{\ell} $ & $ \nu_\ell $ &  $ d_{ORD} $\\
			\hline \hline

			3968&3967&0&1&1&3966&22&2&2&3965&28&2&2\\ \hline
			3968&3964&29&2&2&3963&30&2&2&3962&31&2&2\\ \hline
			3968&3961&32&2&2&3960&33&2&2&3959&44&3&3\\ \hline
			3968&3958&50&4&3&3957&51&4&3&3956&52&4&3\\ \hline
			3968&3955&53&4&3&3954&54&4&3&3953&55&4&3\\ \hline
			3968&3952&56&3&3&3951&57&4&4&3950&58&5&4\\ \hline
			3968&3949&59&6&4&3948&60&7&4&3947&61&8&4\\ \hline
			3968&3946&62&7&4&3945&63&6&4&3944&64&5&4\\ \hline
			3968&3943&65&4&4&3942&66&5&5&3941&72&6&6\\ \hline
			3968&3940&73&6&6&3939&74&6&6&3938&75&6&6\\ \hline
			3968&3937&76&6&6&3936&77&6&6&3935&78&6&6\\ \hline
			3968&3934&79&8&8&3933&80&10&8&3932&81&12&8\\   \hline 
			3968&3931&82&14&8&3930&83&16&8&3929&84&16&8\\   \hline 
			3968&3928&85&16&8&3927&86&16&8&3926&87&16&8\\   \hline 	
			3968&3925&88&17&8&3924&89&14&8&3923&90&14&8\\   \hline 
			3968&3922&91&14&8&3921&92&14&8&3920&93&14&8\\   \hline 
			3968&3919&94&18&8&3918&95&16&8&3917&96&14&8\\   \hline 	
			3968&3916&97&12&8&3915&98&10&8&3914&99&8&8\\   \hline 
			3968&3913&100&9&9& 3912&101&12&12&3911&102&15&15\\   \hline 
			3968&3910&103&18&18&3909&104&21&21&3908&105&24&24\\   \hline 
			3968&3907&106&25&25&3906&107&26&26&3905&108&27&27\\   \hline 
			3968&3904&109&28&28&3903&110&29&28&3902&111&28&28\\   \hline 
			3968&3901&112&29&29&3900&113&30&30&3899&114&31&30\\   \hline 
			3968&3898&115&32&30&3897&116&35&30&3896&117&34&30\\ \hline
			3968&3895&118&33&30&3894&119&32&30&3893&120&31&30\\ \hline
			3968&3892&121&30&30&3891&122&31&31&3890&123&32&32\\ \hline
			3968&3889&124&33&33&3888&125&34&34&3887&126&35&35\\ \hline
			3968&3886&127&36&36&3885&128&37&37&3884&129&38&38\\ \hline
			3968&3883&130&39&39&3882&131&40&40&3881&132&41&41\\ \hline
			3968&3880&133&42&42&3879&134&44&44&3878&135&46&46\\ \hline
			3968&3877&136&48&48&3876&137&50&50&3875&138&52&52\\ \hline
			3968&3874&139&52&52&3873&140&52&52&3872&141&52&52\\ \hline
			3968&3871&142&52&52&3870&143&52&52&3869&144&53&53\\ \hline
			3968&3868&145&54&54&3867&146&55&55&3866&147&56&56\\ \hline
			3968&3865&148&57&57&3864&149&58&58&3863&150&59&59\\   \hline 
			3968&3862&151&60&60&3861&152&61&61&3860&153&62&62\\   \hline 
			3968&3859&154&63&63&3858&155&64&64&3857&156&65&65\\   \hline 
		\end{tabular} 
		
		\caption{Parameters of $C^2_\ell$.}
		\label{tab5}
	\end{table}
	
	\begin{table}[H]
		\centering
		\begin{normalsize}
			\begin{tabular}{|c||c|c|c|c||c|c|c|c||c|c|c|c|}
				\hline
				$ n $ & $ k $ & $ \varrho_{\ell} $ & $ \nu_\ell $ &  $ d_{ORD} $ & $ k $ & $ \varrho_{\ell} $ & $ \nu_\ell $ &  $ d_{ORD} $ & $ k $ & $ \varrho_{\ell} $ & $ \nu_\ell $ &  $ d_{ORD} $\\
				\hline \hline
				3968&3856&157&66&66&3855&158&67&67&3854&159&68&68\\   \hline 
				3968&3853&160&69&69&3852&161&70&70&3851&162&71&71\\   \hline 
				3968&3850&163&72&72&3849&164&73&73&3848&165&74&74\\   \hline 
				3968&3847&166&75&75&3846&167&76&76&3845&168&77&77\\   \hline 
				3968&3844&169&78&78&3843&170&79&79&3842&171&80&80\\   \hline 
				3968&3841&172&81&81&3840&173&82&82&3839&174&83&83\\   \hline 
				3968&3838&175&84&84&3837&176&85&85&3836&177&86&86\\   \hline 
				3968&3835&178&87&87&3834&179&88&88&3833&180&89&89\\   \hline 
				3968&3832&181&90&90&3831&182&91&91&3830&183&92&92\\   \hline 
			\end{tabular} 
		\end{normalsize}
		\caption{Parameters of $C^2_\ell$.}
	\end{table}

	\begin{remark}\label{rem6.3}
		In \cite[Table 1]{GGS} dual AG codes from the first generalized GK curves $\GK_{1,n}$ are constructed. As already recalled, the curves $\GK_{1,n}$ and $\GK_{2,n}$ have the same genus. This allows us to compare the codes obtained in this section with the ones constructed from the curves $\GK_{1,n}$. In particular some of the codes $C^2_\ell$ are shown to have better parameters. The following tables collects some comparisons of our codes from the curve $\GK_{2,n}$ and the ones from the curves $\GK_{1,n}$
	\end{remark}
	\begin{table}[H]
		\centering
		\begin{normalsize}
			\begin{tabular}{|c|c||c|c|}
				\hline
				$ n $ & $ k $ & $ d_{ORD} (\GK_{1,n}) $ & $ d_{ORD} (\GK_{2,n}) $ \\
				\hline \hline
				3968&3910&16&18\\  \hline 
				3968&3909&16&21\\   \hline 
				3968&3908&16&24\\   \hline 
				3968&3907&16&25\\   \hline 
				3968&3906&22&26\\   \hline 
				3968&3905&22&27\\   \hline 
				3968&3904&22&28\\   \hline 
				3968&3903&22&28\\   \hline 
				3968&3902&22&28\\   \hline 
				3968&3901&22&29\\   \hline 
				3968&3900&24&30\\   \hline 
				3968&3899&24&30\\   \hline 
			\end{tabular} 
			\caption{AG codes from the first and second generalized GK curves}
		\end{normalsize}
	\end{table}

	\subsection{AG quantum codes for the second generalized GK curve}
	
	In this section we construct quantum codes  from the curves $\GK_{2,n}$ as an application of the so called \textit{CSS construction} to families of one point AG codes from the curves $\GK_{2,n}$. For more details on quantum codes, we refer the reader to \cite[Section 2]{Quantum}.
	Let $q$ be a prime power. A $q$-ary quantum code of length $N$ and dimension $k$ is defined to be an Hilbert subspace $Q$, with $\dim Q =q^k$, of a $q^n$-dimensional Hilbert space $\mathbb{H}:=(\mathbb{C}^q)^{\otimes n}=\mathbb{C}^q \otimes \cdots \otimes \mathbb{C}^q$. If $Q$ has minimum distance $D$, we will write $Q=[[N,k,D]]_{q^{2n}}$-code.
	
	\begin{proposition} \cite[Lemma 2.5]{Quantum} \label{CSS}
		Let $ C_1 $ and $ C_2 $ be two linear codes with parameters	$ [N, k_i, d_i]_q $, $ i = 1, 2 $, and assume that $ C_1 \subset C_2$. Then there exists an $ [[N, k_2-k_1,D]]_q $-code with $ D = \min\{wt (c) | c \in (C_2 \setminus C_1) \cup (C_1^\perp \setminus C_2^\perp)\} $, where $wt (c)$ is the weight of $c$.
	\end{proposition}
	
	The construction given in Proposition \ref{CSS} is known as \textit{CSS construction}. An application can be obtained looking at the dual of the one point codes from the curves $\GK_{2,n}$. Let $P \in \GK_{2,n}$. Consider $C_2:=C_\ell=C(D,G_2)$ and $C_1:= C_{\ell+s}=C(D,G_1)$, where $s \geq 1$, $G_1=\rho_\ell P$, $G_2=\rho_{\ell+s} P$ and $D=\sum_{Q \in \GK_{2,n}(\mathbb{F}_{q^{2n}}) \setminus \{P\}} Q$. Then we have $C_1 \subset C_2$ and the dimensions of $C_1$ and $C_2$ are  $k_2=N-h_\ell$ and $k_1=N-h_\ell-s$ respectively where $h_i$ denotes the number of non-gaps at $P$ which are smaller than or equal to $i$ and $N= |\GK_{2,n}(\mathbb{F}_{q^{2n}})|-1$. Hence $ k_1-k_2=s$ . From Proposition \ref{CSS} this induces an $[[N,s,D]]_{q^{2n}}$-quantum code, where $D= \min\{wt (c) | c \in (C_2 \setminus C_1) \bigcup (C_1^\perp \setminus C_2^\perp)\}=\min\{wt(c) | c \in (C_\ell \setminus C_{\ell+s})\bigcup (C(P_1+\cdots+P_{N-1},G_1)\setminus C(P_1+\cdots+P_{N-1},G_2))\}$. In particular we get \[D \geq \min \{d_{ORD}(C_\ell),d_1\},\] where $d_1$ denotes the minimum distance of the code $C(D,G_1)$.
	
	Hence the following result follows as a corollary of Proposition \ref{limitcase}.
	
	\begin{corollary}\label{quantum1}
		Let $g=(q-1)(q^{n+1}+q^n-q^2)/2$ and $N=q^{2n+2}-q^{n+3}+q^{n+2}$.
		For every $\ell\in\left[3g-1,N-g\right]$ and $s\in\left[1,N-2\ell\right]$, there exists a quantum code with parameters $[[N,s,D]]_{q^{2n}}$, where $D \geq \ell+1-g$.
		
	\end{corollary}
	
	\begin{proof}
		
		Since there are exactly $g+1$ non-gaps which are less than or equal to $2g$ and $\ell \geq 3g-1$, we have $\rho_{\ell+s}=2g+(\ell+s-(g+1))=g-1+\ell+s$, and hence $d_1 \geq N-\deg(G_1)=N-\rho_{\ell+s}=N-\ell-s-g+1$. Also from $\ell\geq3g-1$ we can apply  Proposition \ref{limitcase} and Proposition \ref{CSS} to get $D \geq \min \{d_{ORD}(C_\ell), d_1\}=\ell+1-g$. The claim follows.
	\end{proof}
	
	Using the tables in Subsection \ref{tables} and applying the general strategy written before, AG quantum codes for which Proposition \ref{limitcase} cannot be applied can also be constructed for $q=2$ and $n=5$.
	Indeed assume in general that $\ell \in [g,3g-1]$. For $s \in [\max\{2g-\ell,1\},N-2\ell]$ we have that $\ell+s \geq 2g$ and $d_1 \geq N-\ell-s-g+1$ as in the proof of Corollary \ref{quantum1}. If $d_{ORD}(C_\ell) \leq  N-\ell-s-g+1$ then arguing as in Corollary \ref{quantum1} there exists a quantum code with parameters $[[N,s,D]]_{q^{2n}}$ where $D \geq d_{ORD}(C_\ell)$. This shows the following proposition.

	\begin{proposition} \label{quantum2}
		Let $\ell \in [g,3g-1]$ and $s \in [\max\{2g-\ell,1\},\min\{N-2\ell,N-\ell-g+1-d_{ORD}(C_\ell)\}]$. Then there exists a quantum code with parameters $[[N,s,D]]_{q^{2n}}$ where $D \geq d_{ORD}(C_\ell)$.
	\end{proposition}
	
	Using the data collected in the tables of Subsection \ref{tables}, the parameters of the AG quantum codes constructed as in Proposition \ref{quantum2} can be determined for $q=2$ and $n=5$. Tables \ref{tab1quantum} and \ref{tab2quantum} collect some explicit examples for $P \in O_1$ and $P \in O_2$ respectively.
	
	\begin{table}[H]
		\centering
		\begin{tabular}{|c|c|c||c|c|c||c|c|c|}
			\hline
			$\ell$ & $d_{Ord}(C_\ell)$ & $s \in$ & $\ell$ & $d_{Ord}(C_\ell)$ & $s \in$ & $\ell$ & $d_{Ord}(C_\ell)$ & $s \in$ \\
			\hline \hline
			
			46&6&$[47,\ldots,3871]$&53&8&$[39,\ldots,3862]$&55&12&$[37,\ldots,3856]$\\ \hline
			57&16& $[35,\ldots,3850]$&59&20&$[33,\ldots,3844]$ &61&24&$[31,\ldots,3838]$\\ \hline
			63&25& $[28,\ldots,3835]$&64&28&$[28,\ldots,3831]$ &75&30&$[17,\ldots,3818]$\\ \hline
			77&32&$[15,\ldots,3814]$&79&34&$[13,\ldots,3810]$ &81&36&$[11,\ldots,3806]$\\ \hline
			83&38&$[9,\ldots,3802]$ &85&39&$[7,\ldots,3798]$ &86&40&$[6,\ldots,3796]$\\   \hline 
			87&41&$[5,\ldots,3794]$&88&42&$[4,\ldots,3792]$&89&44&$[3,\ldots,3790]$\\   \hline 
			91&46&$[1,\ldots,3786]$&93&48&$[1,\ldots,3782]$&95&50&$[1,\ldots,3778]$\\   \hline 
			97&52&$[1,\ldots,3774]$&99&54&$[1,\ldots, 3770]$&101&56&$[1,\ldots,3766]$\\   \hline 
			103&58&$[1,\ldots, 3762]$&105&60&$[1,\ldots,3758]$ &105+i&60+i&$[1,\ldots,3758-i]$\\   \hline 
		\end{tabular}

		\caption{Some $[[3968,s,D \geq d_{ORD}(C_\ell)]]_{2^{10}}$-codes constructed from Proposition \ref{quantum2}, $P \in O_1$.}
		\label{tab1quantum}
	\end{table}

	\begin{table}[H]
		\centering
		\begin{tabular}{|c|c|c||c|c|c|}
			\hline
			$\ell$ & $d_{Ord}(C_\ell)$ & $s \in$ & $\ell$ & $d_{Ord}(C_\ell)$ & $s \in$ \\
			\hline \hline
			46&8&$[46,\ldots,3869]$&82&36&$[10,\ldots,3804]$\\   \hline 
			55&9&$[37,\ldots,3858]$&83&37&$[9,\ldots, 3802]$\\   \hline 
			56&12&$[36,\ldots,3855]$&84&38&$[8,\ldots,3800]$\\   \hline 
			57&15&$[35,\ldots,3851]$&85&39&$[7,\ldots,3798]$\\   \hline 
			58&18&$[34,\ldots,3847]$&86&40&$[6,\ldots,3796]$\\   \hline 
			59&21&$[33,\ldots,3843]$&87&41&$[5,\ldots,3794]$\\   \hline 
			60&24&$[32,\ldots,3839]$&88&42&$[4,\ldots,3792]$\\   \hline 
			61&25&$[31,\ldots,3837]$&89&44&$[3,\ldots,3790]$\\   \hline 
			62&26&$[30,\ldots, 3835]$&90&46&$[2,\ldots,3787]$\\   \hline 
			63&27&$[29,\ldots,3833]$&91&48&$[1,\ldots,3784]$\\   \hline 
			64&28&$[28,\ldots,3831]$&92&50&$[1,\ldots,3781]$\\   \hline 
			67&29&$[25,\ldots,3827]$&93&52&$[1,\ldots,3778]$\\   \hline 
			68&30&$[24,\ldots,3825]$&99&53&$[1,\ldots,3770]$\\   \hline 
			77&31&$[15,\ldots,3814]$&100&54&$[1,\ldots,3768]$\\   \hline 
			78&32&$[14,\ldots,3812]$&101&55&$[1,\ldots,3766]$\\   \hline 
			79&33&$[13,\ldots,3810]$&102&56&$[1,\ldots,3764]$\\   \hline 
			80&34&$[12,\ldots,3808]$&103&57&$[ 1,\ldots,3762]$\\   \hline 
			81&35&$[11,\ldots,3806]$&104&58&$[1,\ldots,3760]$\\   \hline 
		\end{tabular} 
		\caption{Some $[[3968,s,D \geq d_{ORD}(C_\ell)]]_{2^{10}}$-codes constructed from Proposition \ref{quantum2}, $P \in O_2$.}
		\label{tab2quantum}
	\end{table}

	\section*{Acknowledgments}
	
	This research was partially supported by Ministry for Education, University and Research of Italy (MIUR) (Project PRIN 2012 ``Geometrie di Galois e strutture di incidenza'' - Prot. N. 2012XZE22K$_-$005)  and by the Italian National Group for Algebraic and Geometric Structures and their Applications (GNSAGA - INdAM).

	\vspace{1ex}
	\noindent
	Maria Montanucci
	
	\vspace{.5ex}
	\noindent
	Universit\`a degli Studi della Basilicata,\\
	Dipartimento di Matematica Informatica ed Economia,\\
	Contrada Macchia Romana, 85100 Potenza, Italy,\\
	mariamontanucci@gmail.com
	
	\vspace{1ex}
	\noindent
	Vincenzo Pallozzi Lavorante
	
	\vspace{.5ex}
	\noindent
	Universit\'a degli Studi di Modena e Reggio Emilia,\\
	Dipartimento di Matematica Pura e Applicata,\\
	Via Giuseppe Campi 213/b, 41125 Modena, Italy,\\
	vincenzo.pallozzilavorante@unimore.it\\

\end{document}